\documentclass{amsart}

\usepackage[T1]{fontenc}
\usepackage{lmodern}
\usepackage[utf8]{inputenc}
\usepackage[foot]{amsaddr}
\usepackage{amsmath,amssymb,amsthm,mathtools}
\usepackage{geometry}
\usepackage{microtype}
\usepackage{xcolor}
\usepackage{hyperref}
\usepackage{tikz}

\hypersetup{
  colorlinks=true,
  linkcolor=blue,
  citecolor=blue,
  urlcolor=blue
}

\newtheorem{thm}{Theorem}[section]
\newtheorem{cor}[thm]{Corollary}
\newtheorem{prp}[thm]{Proposition}
\newtheorem{lem}[thm]{Lemma}
\theoremstyle{definition}
\newtheorem{dfn}[thm]{Definition}

\DeclareMathOperator{\arccosh}{arccosh}
\DeclareMathOperator{\Aut}{Aut}
\DeclareMathOperator{\diam}{diam}

\newcommand{\B}{\mathcal{B}}
\newcommand{\bbH}{\mathbb{H}}
\newcommand{\calL}{\mathcal{L}}
\newcommand{\D}{\Delta}
\newcommand{\F}{\mathcal{F}}
\newcommand{\op}{\mathrm{op}}
\newcommand{\R}{\mathbb{R}}
\newcommand{\T}{\mathcal{T}}
\newcommand{\Z}{\mathbb{Z}}

\renewcommand{\bar}{\overline}
\renewcommand{\epsilon}{\varepsilon}
\renewcommand{\phi}{\varphi}
\renewcommand{\tilde}{\widetilde}

\title{The Bolza Surface Does Not Locally Minimize the Diameter Function}

\author{Ryo Matsuda}
\address{Ritsumeikan University, Biwako-Kusatsu Campus, 1-1-1 Noji-higashi, Kusatsu, Shiga 525-8577, Japan}
\email{r-mat 'at' fc.ritsumei.ac.jp}

\author{Keita Yamazaki}
\address{Graduate School of Information Sciences, Tohoku University, Aobayama Campus, 6-3-09 Aoba, Aramaki-aza Aoba-ku, Sendai, Miyagi 980-8579, Japan}
\email{yamazaki.keita.s7 'at' dc.tohoku.ac.jp}

\begin{document}

\begin{abstract}
    We study the diameter function on the Teichm\"uller space of closed hyperbolic Riemann surfaces.
    Although the Bolza surface maximizes the systole among hyperbolic surfaces of genus \(2\),
    its extremal behavior with respect to the diameter function is less clear.
    Using first variations of distances at diameter pairs and linear programming duality,
    we establish a sufficient criterion for a surface not to be a local minimizer of the diameter function.
    More precisely,
    we show that if each diameter pair can be assigned a simple closed geodesic formed by two minimizing geodesic segments joining its endpoints,
    with distinct assigned geodesics pairwise disjoint,
    then a deformation varying only the length parameters in Fenchel--Nielsen coordinates strictly decreases the diameter.
    The assigned geodesics extend to a pants decomposition,
    and the decrease is detected by a strictly negative upper bound for the first-order variation of the diameter.
    We then apply this criterion to the Bolza surface and prove that it is not a local minimizer of the diameter function on the Teichm\"uller space of genus \(2\).
\end{abstract}

\maketitle

\section{Introduction}

Let \(S_g\) be a closed oriented smooth surface of genus \(g \ge 2\).
A marked Riemann surface is a pair \((X, \phi)\),
where \(X\) is a closed Riemann surface of genus \(g\) and \(\phi \colon S_g \to X\) is an orientation-preserving homeomorphism.
Two marked Riemann surfaces \((X_1, \phi_1)\) and \((X_2, \phi_2)\) are said to be marking equivalent if there exists a biholomorphic map \(h \colon X_1 \to X_2\) such that \(h \circ \phi_1\) is homotopic to \(\phi_2\).
The Teichm\"uller space \(\T_g\) of \(S_g\) is the set of marking equivalence classes of marked Riemann surfaces.
We often write \(X\) for \([X, \phi]\) when the marking is understood.

Every Riemann surface representing a point of \(\T_g\) admits a unique hyperbolic metric of curvature \(- 1\) in its conformal class.
A pants decomposition \(\Gamma = \{\gamma_k\}_{k = 1}^{3 g - 3}\) of \(S_g\) determines Fenchel--Nielsen coordinates \((l_k, \tau_k)_{k = 1}^{3 g - 3}\) on \(\T_g\),
where \(l_k\) and \(\tau_k\) denote the length and twist parameters,
respectively (see \cite{Buser}).
The length parameters can be varied independently while keeping the twist parameters fixed.
For an essential simple closed curve \(\gamma\) on \(S_g\),
let \(l_\gamma \colon \T_g \to \R_{> 0}\) denote its hyperbolic length function.
Let \(d_X(p, q)\) denote the hyperbolic distance between two points \(p\)
and \(q\) on \(X\).
The diameter of \(X\) is the conformal invariant defined by
\[
    \diam(X) \coloneqq \max_{p, q \in X} d_X(p, q).
\]
Thus the diameter can be regarded as a function \(\diam \colon \T_g \to \R_{> 0}\) on Teichm\"uller space.
The problem of minimizing this function asks how small the greatest distance between two points can be made among all surfaces of a fixed genus.

Stepanyants--Beardon--Paton--Krioukov \cite{SBPK} computed the diameters of generalized Bolza surfaces exactly.
In particular,
they showed that the Bolza surface \(B\) satisfies
\[
    \diam(B) = \arccosh(3 + 2 \sqrt{2}).
\]
Irmer--Saha \cite{IrmerSaha} proved that the diameter function is a topological Morse function on Teichm\"uller space.
These results provide a framework for studying the local behavior of the diameter function and the geometry of its extrema.
Among Riemann surfaces of genus \(2\),
the Bolza surface has an automorphism group of maximal order,
and it is known to be a maximizer of the systole function,
which assigns to a surface the length of its shortest closed geodesic,
on \(\T_2\) \cite{SchmutzSchaller1999}.
By a result of Chang \cite{Chang1977},
increasing the systole decreases an upper bound for the diameter when the genus is fixed.
These observations raise the question of whether the Bolza surface minimizes the diameter function on \(\T_2\).

The purpose of this paper is to establish a geometric criterion for the existence of a deformation that decreases the diameter and to apply this criterion to the Bolza surface.
Our main results are as follows.

\begin{thm}\label{thm:main}
    Let \(X \in \T_g\) be represented by a marked closed hyperbolic Riemann surface of genus \(g \ge 2\).
    Define its diameter locus by
    \[
        \calL(X) \coloneqq \{\{p, q\} \subset X \mid d_X(p, q) = \diam(X)\}.
    \]
    Suppose that for each \(P \in \calL(X)\) one can choose two minimizing geodesic segments \(a_P^+\) and \(a_P^-\) joining the points of \(P\) so that \(\gamma_P \coloneqq a_P^+ \cup a_P^-\) is a simple closed geodesic and,
    whenever \(P \neq Q\),
    either \(\gamma_P = \gamma_Q\) or \(\gamma_P \cap \gamma_Q = \emptyset\).
    Then the set of distinct geodesics \(\Gamma_0 \coloneqq \{\gamma_P \mid P \in \calL(X)\}\) can be extended to a pants decomposition of \(X\).
    With respect to the associated Fenchel--Nielsen coordinates,
    there exists a smooth curve \(t \mapsto X_t \in \T_g\),
    defined for real \(t\) near \(0\),
    with \(X_0 = X\),
    that decreases the length parameters corresponding to the geodesics in \(\Gamma_0\),
    while keeping all other Fenchel--Nielsen coordinates fixed.
    Writing \(V = \dot{X}_0 \in T_X \T_g\),
    we have
    \[
        \lim_{t \downarrow 0} \frac{\diam(X_t) - \diam(X)}{t} \le \frac{1}{2} \max_{\gamma \in \Gamma_0} dl_\gamma(V) < 0,
    \]
    and,
    for sufficiently small \(t > 0\),
    \[
        \diam(X_t) < \diam(X).
    \]
    Therefore \(X\) is not a local minimizer of the diameter function on \(\T_g\).
\end{thm}

\begin{cor}\label{cor:bolza}
    The Bolza surface is not a local minimizer of the diameter function on \(\T_2\).
\end{cor}

This result suggests that the systole and diameter functions exhibit essentially different behavior as topological Morse functions.
It thus provides further motivation for studying the diameter function in its own right.

\section*{Acknowledgements}
The authors thank Toshiyuki Sugawa and Hidetoshi Masai for helpful discussions,
and Kento Sakai and Ryutaro Misawa for providing useful information.

\section{First Variations at Diameter Pairs}\label{sec:first-variations-at-diameter-pairs}

Throughout this section,
let \(X \in \T_g\) be a closed hyperbolic Riemann surface,
and let \(D = \diam(X)\).
Recall that the diameter locus of \(X\) is the set
\[
    \calL(X) = \{\{p, q\} \subset X \mid d_X(p, q) = D\}.
\]
We call its elements diameter pairs.
Fix a diameter pair \(P = \{p, q\} \in \calL(X)\),
and denote by \(a_1, \ldots, a_m\) all the minimizing geodesic segments joining \(p\) and \(q\),
each oriented from \(p\) to \(q\).
We verify on the universal cover that there are only finitely many minimizing geodesic segments joining \(p\) and \(q\),
and that their common length is separated from the lengths of all other geodesic segments by a positive gap.
Write \(X = \bbH^2 / \F\),
where \(\F\) is a torsion-free cocompact Fuchsian group,
and fix lifts \(\tilde{p}\) and \(\tilde{q}\) of \(p\) and \(q\).
The relative homotopy classes of paths from \(p\) to \(q\) correspond to the geodesic segments joining \(\tilde{p}\) to \(f(\tilde{q})\),
\(f \in \F\).
Since \(\F\) acts properly discontinuously,
for every \(L > 0\) there are only finitely many \(f \in \F\) such that \(d_{\bbH^2}(\tilde{p}, f(\tilde{q})) \le L\).
Hence there are finitely many minimizing segments,
and their common length \(D\) is separated from the lengths of all nonminimizing segments by some \(\delta > 0\).

Next,
fix \(L > D\).
Among the finitely many segments whose lengths lie in \((D, L]\),
take the differences between their lengths and \(D\),
and let \(\delta\) be the minimum of these differences and \(L - D\).
If no such segment exists,
set \(\delta = L - D\).
Then \(\delta > 0\),
and every segment that is not minimizing has length at least \(D + \delta\).

Let \(t \mapsto X_t\) be a \(C^1\) curve in \(\T_g\),
defined for real \(t\) in a neighborhood of \(0\),
with \(X_0 = X\).
Let \(V = \dot{X}_0 \in T_X \T_g\) denote its tangent vector at \(t = 0\).
Choose a smooth local section of the space of hyperbolic metrics over \(\T_g\) near \(X\).
This section represents \(X_t\) by a \(C^1\) family of smooth hyperbolic metrics \(g_t\) on a fixed smooth surface \(S\),
with \((S, g_0) = X\).

Choose disjoint normal neighborhoods \(U_p\) of \(p\) and \(U_q\) of \(q\).
For \((p', q') \in U_p \times U_q\),
identify relative homotopy classes by adjoining the short geodesic arcs at the two endpoints.
For each relative homotopy class \(\alpha\), let \(L_\alpha(t, p', q')\) be the length of its geodesic representative for \(g_t\).
For the relative homotopy classes represented by \(a_1, \ldots, a_m\),
write \(L_j = L_{a_j}\),
\(j = 1, \ldots, m\).
The unique lifted geodesic between distinct endpoints has no conjugate points.
Smooth dependence of the geodesic equation and the inverse function theorem therefore give continuous joint first derivatives of these functions in \((t, p', q')\), and smooth dependence on the endpoints for each fixed \(t\).

We justify uniformly that no initially nonminimizing class can become minimizing nearby.
Fix \(\rho < \delta / 4\) and restrict the endpoints so that
\(d_{g_0}(p, p') + d_{g_0}(q, q') \le \rho\).
Adjoining the short endpoint arcs, and then reversing this construction, gives
\[
    |L_\alpha(0, p', q') - L_\alpha(0, p, q)| \le \rho
\]
for every class \(\alpha\).
Uniform comparison of the metrics on the compact surface gives, for small \(|t|\),
\[
    (1 - \epsilon) L_\alpha(0, p', q') \le L_\alpha(t, p', q') \le (1 + \epsilon) L_\alpha(0, p', q'),
\]
where \(\epsilon > 0\) can be chosen arbitrarily small by restricting \(t\), uniformly in the class and the endpoints.
An initially nonminimizing class thus has length at least \((1 - \epsilon) (D + \delta - \rho)\),
whereas every initially minimizing class has length at most \((1 + \epsilon) (D + \rho)\).
Choose \(\epsilon\) so small that the former bound exceeds the latter.
This proves
\begin{equation}\label{eq:finite-min}
    d_{g_t}(p', q') = \min_{1 \le j \le m} L_j(t, p', q')
\end{equation}
on a fixed neighborhood of \((0, p, q)\).

Let \(u_j \in T_p X\) and \(v_j \in T_q X\) be the initial and terminal unit tangent vectors of \(a_j\),
respectively.
Define the metric variation with fixed endpoints by
\[
    \lambda_j(V) \coloneqq \left.\frac{d}{dt}\right|_{t = 0} L_j(t, p, q).
\]
In the chosen local section,
each \(\lambda_j\) is a linear functional on \(T_X \T_g\).
The functionals \(\lambda_j\) may depend on that section;
the quantity \(\Phi_P\) below does not,
as we shall verify.

If the endpoints move according to
\[
    p' = \exp_p(t \xi), \qquad q' = \exp_q(t \eta),
\]
where \(\xi \in T_p X\) and \(\eta \in T_q X\),
then the first variation formula for length,
allowing both the endpoints and the metric to vary,
gives
\[
  L_j(t, p', q') = D + t (\lambda_j(V) - \langle u_j, \xi\rangle + \langle v_j, \eta \rangle) + o(t).
\]

To study the first variation of the local maximum of the distance near \(P\),
we introduce the following quantity.

\begin{dfn}
    For each \(j = 1, \ldots, m\),
    put
    \[
        b_j \coloneqq (- u_j, v_j) \in T_p X \oplus T_q X.
    \]
    For \(V \in T_X \T_g\),
    define the local directional quantity of \(P\) by
    \begin{equation}\label{eq:Phi}
        \Phi_P(V) = \sup_{\zeta \in T_p X \oplus T_q X} \min_{1 \le j \le m} (\lambda_j(V) + \langle b_j, \zeta \rangle).
    \end{equation}
\end{dfn}

\begin{prp}\label{prop:LP}
    For every diameter pair \(P = \{p, q\}\) and every \(V \in T_X \T_g\),
    we have
    \[
        \Phi_P(V) = \min_{\mu \in \B(P)} \sum_{j = 1}^m \mu_j \lambda_j(V),
    \]
    where
    \[
        \B(P) \coloneqq \left\{\mu \in \R_{\ge 0}^m\; \middle|\; \sum_{j = 1}^m \mu_j = 1,\ \sum_{j = 1}^m \mu_j b_j = 0\right\}.
    \]
\end{prp}

\begin{proof}
    The quantity \(\Phi_P(V)\) defined in \eqref{eq:Phi} is the optimal value of the following finite-dimensional linear program:
    \[
        \begin{array}{ll}
            \text{maximize}   & s,\\
            \text{subject to} & s - \langle b_j, \zeta \rangle \le \lambda_j(V), \qquad j = 1, \ldots, m,
        \end{array}
    \]
    with variables \(s \in \R\) and \(\zeta \in T_p X \oplus T_q X\). 
    Therefore the dual problem is
    \[
        \begin{array}{ll}
            \text{minimize}   & \displaystyle \sum_{j = 1}^m \mu_j \lambda_j(V),\\
            \text{subject to} & \mu_j \ge 0, \qquad j = 1, \ldots, m,\\
                              & \displaystyle \sum_{j = 1}^m \mu_j = 1, \qquad \sum_{j = 1}^m \mu_j b_j = 0.
        \end{array}
    \]
    Thus its feasible set is precisely \(\B(P)\).

    The primal problem is feasible by taking \(\zeta = 0\) and \(s\) sufficiently small.
    For \(\zeta = (\xi, \eta)\), the first variation at fixed metric and \eqref{eq:finite-min} give,
    as \(h \downarrow 0\),
    \[
        d_X(\exp_p(h \xi), \exp_q(h \eta)) = D + h \min_{1 \le j \le m} \langle b_j, \zeta \rangle + o(h).
    \]
    Since the left-hand side is at most \(D\),
    we have \(\min_j \langle b_j, \zeta \rangle \le 0\) for every \(\zeta\).
    Thus every feasible \((s, \zeta)\) satisfies
    \[
        s \le \min_{1 \le j \le m}(\lambda_j(V) + \langle b_j, \zeta \rangle) \le \max_{1 \le j \le m} \lambda_j(V) + \min_{1 \le j \le m} \langle b_j, \zeta \rangle \le \max_{1 \le j \le m} \lambda_j(V).
    \]
    By strong duality for finite-dimensional linear programming \cite[Chapter 6]{MatousekGartner},
    the dual minimum is attained and equals the primal optimal value.
    Hence
    \[
        \Phi_P(V) = \min_{\mu \in \B(P)} \sum_{j = 1}^m \mu_j \lambda_j(V).
    \]
\end{proof}

In particular,
\(\Phi_P(V)\) is finite and \(\B(P)\) is nonempty.
The defining constraints make \(\B(P)\) a closed subset of the standard probability simplex,
so it is compact.

\begin{lem}\label{lem:diameter-pair-estimates}
    The diameter locus \(\calL(X)\) is finite,
    and for each diameter pair \(P = \{p, q\}\),
    there is a constant \(\kappa_P > 0\) such that
    \begin{equation}\label{eq:endpoint-estimate}
        \min_{1 \le j \le m} \langle b_j, \zeta \rangle \le - \kappa_P \|\zeta\|, \qquad \zeta \in T_p X \oplus T_q X,
    \end{equation}
    where the norm is induced by the product metric.
    For all sufficiently small \(\zeta = (\xi, \eta) \in T_p X \oplus T_q X\),
    we have 
    \[
        d_X(\exp_p(\xi), \exp_q(\eta)) \le D - \frac{\kappa_P}{2} \|\zeta\|.
    \]
\end{lem}

\begin{proof}
    Since \(X\) is closed and has sectional curvature \(- 1\),
    the local maxima of \(d_X\) on \(X \times X\) are isolated by \cite[Remark 3.6]{Horja1997}.
    Thus the set \(d_X^{- 1}(D)\) of ordered diameter pairs is discrete.
    It is also closed in the compact space \(X \times X\),
    so it is finite.
    Forgetting the order gives the finiteness of \(\calL(X)\).

    Fix \(P = \{p, q\}\) and write \(F_P(\zeta) = \min_j \langle b_j, \zeta \rangle\).
    Suppose that \(F_P(\zeta) \ge 0\) for some nonzero \(\zeta = (\xi, \eta)\).
    For \(s\) sufficiently close to \(0\),
    set \(p_s = \exp_p(s \xi)\),
    \(q_s = \exp_q(s \eta)\),
    and \(f_j(s) = L_j(0, p_s, q_s)\).
    For each \(j\),
    choose lifts of the two endpoint paths to \(\bbH^2\) whose initial points are joined by a lift of \(a_j\).
    Then \(f_j(s)\) is the distance between these lifted paths at time \(s\).
    Since both paths are geodesics,
    \(f_j\) is convex near \(0\) by \cite[Lemma 2.3]{Horja1997}.
    The endpoint first variation gives \(f_j'(0) = \langle b_j, \zeta \rangle \ge 0\).
    Therefore \(f_j(s) \ge f_j(0) + s f_j'(0) \ge D\) for all sufficiently small \(s > 0\).
    By \eqref{eq:finite-min},
    \(D \ge d_X(p_s, q_s) = \min_j f_j(s) \ge D\).
    Since \(\zeta \neq 0\),
    this contradicts the isolation of \((p, q)\).

    Thus \(F_P(\zeta) < 0\) for every \(\zeta \neq 0\).
    On the compact unit sphere,
    the continuous function \(F_P\) attains a negative maximum,
    which we denote by \(- \kappa_P\).
    Since \(F_P(h \zeta) = h F_P(\zeta)\) for \(h \ge 0\),
    this gives \eqref{eq:endpoint-estimate}.
    Finally,
    the first-order expansions of the finitely many branches in \eqref{eq:finite-min} give
    \[
        d_X(\exp_p(\xi), \exp_q(\eta)) = D + F_P(\zeta) + o(\|\zeta\|), \qquad \zeta = (\xi, \eta) \to 0.
    \]
    The local distance estimate follows.
\end{proof}

To verify that \(\Phi_P(V)\) is independent of the choice of smooth local section,
consider two such sections,
after identifying their reference surfaces so that the corresponding metrics agree at \(t = 0\).
If the endpoints fixed in the new representation have velocity \(\zeta' \in T_p X \oplus T_q X\) with respect to the original representation,
the first variation formula gives
\[
    \tilde{\lambda}_j(V) = \lambda_j(V) + \langle b_j, \zeta' \rangle.
\]
Thus,
by translating the variable \(\zeta\),
\[
    \begin{aligned}
        \tilde{\Phi}_P(V) &= \sup_{\zeta \in T_p X \oplus T_q X} \min_{1 \le j \le m} (\lambda_j(V) + \langle b_j, \zeta + \zeta' \rangle)\\
                          &= \Phi_P(V).
    \end{aligned}
\]
Hence \(\Phi_P\) is independent of the choice of this local representation.

We use the fixed smooth surface \(S\) to compare the metrics \(g_t\).
For \(\zeta = (\xi, \eta)\) near \(0\) in \(T_p X \oplus T_q X\),
define
\[
    d_t(\zeta) \coloneqq d_{g_t}(\exp_p(\xi), \exp_q(\eta)),
\]
where both exponential maps are taken with respect to \(g_0\).

\begin{lem}\label{lem:local-max}
    Let \(P = \{p, q\}\) be a diameter pair.
    Let \(U \subset U_p \times U_q\) be the image of a sufficiently small open ball about \(0\) under the product exponential map at \((p, q)\).
    Put
    \[
        M_P(t) \coloneqq \max_{(p', q') \in \bar{U}} d_{g_t}(p', q').
    \]
    Then,
    for sufficiently small \(|t|\),
    all maximizers lie in the interior of \(U\),
    and 
    \begin{equation}\label{eq:local-asymptotic}
        M_P(t) = D + t \Phi_P(V) + o(t) \qquad (t \downarrow 0).
    \end{equation}
\end{lem}

\begin{proof}
    Use the product exponential coordinates \(\zeta = (\xi, \eta)\) on \(U\).
    Equation \eqref{eq:finite-min} and the uniform first-order expansions of the functions \(L_j\) give
    \begin{equation}\label{eq:uniform-expansion}
        d_t(\zeta) = D + \min_{1 \le j \le m} (t \lambda_j(V) + \langle b_j, \zeta \rangle) + o(|t| + \|\zeta\|).
    \end{equation}
    By Lemma \ref{lem:diameter-pair-estimates},
    with \(\kappa = \kappa_P\),
    inequality \eqref{eq:endpoint-estimate} holds and,
    for all sufficiently small \(\zeta\),
    \[
        d_0(\zeta) \le D - \frac{\kappa}{2} \|\zeta\|.
    \]
    
    After shrinking \(U\),
    choose \(t_0 > 0\) so that all \(L_j\) are defined on a neighborhood of the compact set \([- t_0, t_0] \times \bar{U}\).
    By continuity,
    the partial derivatives \(\partial_t L_j\) are uniformly bounded on this compact set.
    Let \(r > 0\) be the coordinate radius of the resulting neighborhood \(U\).
    Since the minimum of finitely many numbers changes by at most the largest change of an entry,
    there exists \(C_0 > 0\) such that
    \[
        |d_t(\zeta) - d_0(\zeta)| \le C_0 |t| \qquad (\zeta \in \bar{U}).
    \]
    On the boundary,
    \(d_0(\zeta) \le D - \kappa r / 2\),
    whereas \(d_t(0) \ge D - C_0 |t|\).
    Hence,
    for sufficiently small \(|t|\),
    every maximizer on \(\bar{U}\) lies in the interior of \(U\).
    If \(\zeta_t\) is such a maximizer,
    then
    \[
        D - C_0 |t| \le d_t(0) \le d_t(\zeta_t) \le D - \frac{\kappa}{2} \|\zeta_t\| + C_0 |t|.
    \]
    Therefore
    \[
        \|\zeta_t\| \le \frac{4 C_0}{\kappa} |t|.
    \]
    
    Put
    \[
        H_V(\omega) \coloneqq \min_{1 \le j \le m} (\lambda_j(V) + \langle b_j, \omega \rangle).
    \]
    Inequality \eqref{eq:endpoint-estimate} gives
    \[
        H_V(\omega) \le \max_{1 \le j \le m} \lambda_j(V) - \kappa \|\omega\|.
    \]
    Thus \(H_V\) tends to \(- \infty\) at infinity,
    and \(\Phi_P(V) = \sup_\omega H_V(\omega)\) is attained at some \(\omega_*\).
    
    For \(t > 0\),
    setting \(\zeta = t \omega\) in \eqref{eq:uniform-expansion} gives
    \[
        d_t(t \omega) = D + t \min_{1 \le j \le m} (\lambda_j(V) + \langle b_j, \omega \rangle ) + o(t),
    \]
    where the \(o(t)\) term is uniform for \(\omega\) in any fixed bounded subset of \(T_p X \oplus T_q X\).
    Since \(\omega_t \coloneqq \zeta_t / t\) is bounded,
    \[
        M_P(t) = D + t H_V(\omega_t) + o(t) \le D + t \Phi_P(V) + o(t).
    \]
    On the other hand,
    \(t \omega_* \in U\) for sufficiently small \(t > 0\),
    so
    \[
        M_P(t) \ge d_t(t \omega_*) = D + t \Phi_P(V) + o(t).
    \]
    Combining the two bounds yields \eqref{eq:local-asymptotic}.
\end{proof}

\section{Fenchel--Nielsen Deformations}

In this section,
we prove Theorem \ref{thm:main}.

\begin{prp}\label{prp:certificate}
    Let \(P = \{p, q\}\) be a diameter pair,
    and let \(a^+\) and \(a^-\) be two minimizing geodesic segments from \(p\) to \(q\).
    If \(\gamma \coloneqq a^+ \cup a^-\) is a simple closed geodesic,
    then,
    for every \(V \in T_X \T_g\),
    \begin{equation}\label{eq:certificate}
        \Phi_P(V) \le \frac{1}{2} dl_\gamma(V).
    \end{equation}
\end{prp}

\begin{proof}
    Let \(b^+, b^- \in T_p X \oplus T_q X\) be the vectors associated with \(a^+\) and \(a^-\),
    respectively,
    as defined in Section \ref{sec:first-variations-at-diameter-pairs}.
    Since \(\gamma = a^+ \cup a^-\) is a simple closed geodesic,
    \(a^+\) and \(a^-\) are its two complementary subarcs between \(p\) and \(q\).
    When both segments are oriented from \(p\) to \(q\),
    their tangent directions are opposite at both endpoints.
    Hence
    \[
        b^- = - b^+.
    \]
    
    Therefore the measure assigning weight \(1 / 2\) to each of these two segments and weight \(0\) to all the other minimizing segments belongs to \(\B(P)\).
    Proposition \ref{prop:LP} gives
    \[
        \Phi_P(V) \le \frac{1}{2} (\lambda^+(V) + \lambda^-(V)).
    \]
    For the family of metrics \(g_t\) on the fixed surface,
    write \(\dot{g}_0 \coloneqq \left.\frac{d}{dt}\right|_{t = 0} g_t\),
    and let \(T\) be a \(g_0\)-unit tangent vector field along \(\gamma\),
    with \(ds\) denoting the \(g_0\)-length element.
    The first variation formula with fixed endpoints gives
    \[
        \lambda^\pm(V) = \frac{1}{2} \int_{a^\pm} \dot{g}_0(T, T) \, ds.
    \]
    Since the initial curve \(\gamma\) is a closed geodesic,
    the first variation formula for closed geodesic length gives
    \begin{equation}
        \lambda^+(V) + \lambda^-(V) = \frac{1}{2} \int_\gamma \dot{g}_0(T, T) \, ds = dl_\gamma(V).
    \end{equation}
    This proves \eqref{eq:certificate}.
\end{proof}

\begin{proof}[Proof of Theorem \ref{thm:main}]
    Put \(D = \diam(X)\).
    By Lemma \ref{lem:diameter-pair-estimates},
    we may write
    \[
        \calL(X) = \{P_1, \ldots, P_N\}.
    \]
    Write \(P_i = \{p_i, q_i\}\) and put \(\Gamma_0 = \{\gamma_{P_i}\}_{1 \le i \le N}\).
    By hypothesis,
    the distinct geodesics in \(\Gamma_0\) are pairwise disjoint.
    Since distinct simple closed geodesics on \(X\) represent distinct essential free homotopy classes,
    \(\Gamma_0\) can be extended to a pants decomposition \(\Gamma\) of \(X\).

    Let \((l_\gamma, \tau_\gamma)_{\gamma \in \Gamma}\) be the Fenchel--Nielsen coordinates associated with this pants decomposition.
    For sufficiently small real \(t\),
    define a smooth curve \(X_t\) with \(X_0 = X\) by
    \[
        l_\gamma(X_t) = \begin{cases}
                            l_\gamma(X) - t, & \gamma \in \Gamma_0,\\
                            l_\gamma(X),     & \gamma \in \Gamma \setminus \Gamma_0,
                        \end{cases}
        \qquad
        \tau_\gamma(X_t) = \tau_\gamma(X), \quad \gamma \in \Gamma.
    \]
    Writing \(V = \dot{X}_0\),
    we have
    \[
        dl_{\gamma_{P_i}}(V) = - 1, \qquad i = 1, \ldots, N.
    \]
    Applying Proposition \ref{prp:certificate} to each diameter pair \(P_i\) yields
    \begin{equation}\label{eq:phi-negative}
        \Phi_{P_i}(V) \le \frac{1}{2} dl_{\gamma_{P_i}}(V) = - \frac{1}{2}.
    \end{equation}

    Represent \(X_t\) by metrics \(g_t\) on the fixed surface \(S\).
    For each \(i\),
    choose a sufficiently small coordinate ball \(W_i\) about \((p_i, q_i)\) in \(S \times S\) as in Lemma \ref{lem:local-max}.
    Put
    \[
        M_i(t) = \max_{(x, y) \in \bar{W_i}} d_{g_t}(x, y).
    \]
    By Lemma \ref{lem:local-max},
    \begin{equation}\label{eq:general-local-max}
        M_i(t) = D + t \Phi_{P_i}(V) + o(t) \qquad (t \downarrow 0).
    \end{equation}
    Let \(W_i^\op = \{(y, x) \mid (x, y) \in W_i\}\).
    By symmetry of the distance,
    the maximum on \(\bar{W_i^\op}\) is also \(M_i(t)\).
    
    The compact set
    \[
        K = (S \times S) \setminus \bigcup_{i = 1}^N(W_i \cup W_i^\op)
    \]
    contains no ordered pair realizing \(D\),
    so \(d_{g_0} \le D - \epsilon\) on \(K\) for some \(\epsilon > 0\).
    Since \(g_t \to g_0\) uniformly on \(S\),
    the distance functions converge uniformly on \(S \times S\).
    Thus \(d_{g_t} \le D - \epsilon / 2\) on \(K\) for all sufficiently small \(|t|\).
    Since \(M_i(t) \to D\) as \(t \downarrow 0\),
    it follows that,
    for all sufficiently small \(t > 0\),
    \[
        \diam(X_t) = \max_{1 \le i \le N} M_i(t).
    \]
    
    Hence \eqref{eq:phi-negative} and \eqref{eq:general-local-max} imply
    \begin{equation}
        \begin{aligned}
            \lim_{t \downarrow 0} \frac{\diam(X_t) - D}{t} & = \max_{1 \le i \le N} \Phi_{P_i}(V)\\
                                                           & \le \frac{1}{2} \max_{\gamma \in \Gamma_0} dl_\gamma(V)\\
                                                           & = - \frac{1}{2} < 0.
        \end{aligned}
    \end{equation}
    Since \(N\) is finite,
    the remainders in the local asymptotic expansions can be chosen uniformly in \(i\),
    which justifies the first equality.
\end{proof}

\section{Application to the Bolza Surface}

Let \(B\) denote the Bolza surface obtained by identifying opposite sides of a regular hyperbolic octagon \(P\) with interior angle \(\pi / 4\).
Let \(s \coloneqq 2 \arccosh(1 + \sqrt{2})\) be its side length.
Write \(O\) for the image on the Bolza surface of the origin,
the center of this fundamental domain,
and \(V\) for the common image its vertices \(v_0, \ldots, v_7\),
respectively.
If \(R\) is the circumradius of this fundamental domain,
then,
by \cite[Theorem 1]{SBPK},
\[
    D_B \coloneqq \diam(B) = R = \arccosh(3 + 2 \sqrt{2}).
\]

The points \(O\) and \(V\) are Weierstrass points,
as are the images of the side midpoints
\[
    \tanh \frac{s}{4} e^{j \pi i/ 4}, \qquad j = 0, 1, 2, 3,
\]
which we denote by \(M_j\).
The following description is stated in \cite[Example 6]{IrmerSaha}.
We give an elementary proof in the unit disk using the distance estimate of \cite{SBPK}.

\begin{prp}\label{prp:bolza-pairs}
    The diameter locus \(\calL(B)\) of the Bolza surface \(B\) is
    \[
        \calL(B) = \{\{O, V\}, \{M_0, M_2\}, \{M_1, M_3\}\}.
    \]
\end{prp}

\begin{proof}
    Work in the unit disk \(\D\),
    write \(\delta\) for its hyperbolic distance and \(d\) for the distance on \(B\),
    and let \(\pi \colon \D \to B\) be the universal covering map.
    The regular octagon \(P\),
    centered at \(0\),
    has inradius \(s / 2\) and circumradius \(R\),
    with
    \[
        \cosh \frac{s}{2} = 1 + \sqrt{2}, \qquad \cosh R = \cosh^2 \frac{s}{2}, \qquad R < s.
    \]

    We first recall a geometric estimate from \cite{SBPK} that we use.
    For adjacent vertices \(v_1\),
    \(v_2\),
    let \(T = \triangle 0 v_1 v_2\) and let \(Q_i\) be the dual octagon centered at \(v_i\).
    The octagons \(P\) and \(Q_i\) are Dirichlet domains for \(O\) and \(V\),
    respectively;
    a dual symmetry interchanges \(O\) and \(V\).
    As in \cite[Proof of Proposition 1]{SBPK},
    there is a deck transformation \(h\) translating along \(v_1 v_2\) by a distance \(s\),
    with \(h(v_1) = v_2\) and \(h(Q_1) = Q_2\).
    Let \(w_1 \in Q_1\) satisfy \(b \coloneqq \delta(v_1, w_1) \le s / 2\),
    and put \(w_2 \coloneqq h(w_1) \in Q_2\).
    Then \(\delta(v_2, w_2) = b\),
    and
    \begin{equation}\label{eq:paired-distance-bound}
        z \in T,\ \delta(0, z) \ge b \qquad \implies \qquad f(z) \coloneqq \min\{\delta(z, w_1), \delta(z, w_2)\} \le R.
    \end{equation}
    These facts follow from \cite[Proposition 1 and Theorems 3 and 4]{SBPK}.
    Boundary points are included by continuity.

    Let \(\alpha(t)\) and \(\beta(t)\) be unit-speed geodesics with \(\delta(\alpha(0), \beta(0)) = R\).
    By applying an isometry of \(\D\),
    we may assume that \(\alpha(0) = 0\) and \(\beta(0) = a = \tanh(R / 2)\).
    With \(u = \tanh(t / 2)\),
    write
    \[
        \alpha(t) = u e^{i \theta}, \qquad \beta(t) = \frac{a + u e^{i \psi}}{1 + a u e^{i \psi}}.
    \]
    A direct calculation in the disk model gives,
    for \(t > 0\),
    \begin{equation}\label{eq:distance-difference-identity}
        \frac{\cosh \delta(\alpha(t), \beta(t)) - \cosh R}{\sinh t \cosh t} = H \tanh t + K,
    \end{equation}
    where
    \begin{align*}
        &H \coloneqq (\cosh R - 1)(1 - \cos \theta \cos \psi) + 1 - \cos(\theta - \psi) \ge 0,\\
        &K \coloneqq \sinh R (\cos \psi - \cos \theta).
    \end{align*}
    Thus the right-hand side is nondecreasing in \(t\).
    It is strictly increasing unless \(H = 0\),
    which occurs precisely when both geodesics move along the same line in the same direction;
    in that case \(K = 0\) and their distance is constantly \(R\).
    
    For \(0 < b < s / 2\),
    fix \(w_1\),
    \(w_2\) and put \(\Omega_b \coloneqq \{z \in T \mid \delta(0, z) \ge b\}\).
    The maximum of \(f\) on \(\Omega_b\) occurs on its boundary \cite[Theorem 3]{SBPK}.

    On a radial side \(0 v_i\),
    if \(\delta(0, z) > b\),
    the triangle inequality gives
    \[
        f(z) \le \delta(z, w_i) \le R - \delta(0, z) + b < R.
    \]
    
    On the side \(v_1 v_2\),
    put \(x = \delta(v_1, z)\).
    For \(0 < x < s\),
    let \(\theta\) be the angle between \(v_1 w_1\) and \(v_1 v_2\).
    Since \(w_2 = h(w_1)\),
    the hyperbolic cosine rule and a weighted average estimate give
    \[
        \begin{aligned}
            \cosh \delta(z, w_1) &= \cosh b \cosh x - \sinh b \sinh x \cos \theta,\\
            \cosh \delta(z, w_2) &= \cosh b \cosh(s - x) + \sinh b \sinh(s - x) \cos \theta,\\
            \min_{i = 1, 2} \cosh \delta(z, w_i) &\le \frac{\sinh(s - x) \cosh \delta(z, w_1) + \sinh x \cosh \delta(z, w_2)}{\sinh(s - x) + \sinh x}\\
                                                 &= \frac{\cosh b \cosh(s / 2)}{\cosh(x - s / 2)} < \cosh R.
        \end{aligned}
    \]
    At either vertex,
    \(f(z) \le b < R\).

    It remains to consider \(\delta(0, z) = b\).
    Suppose that both distances are at least \(R\).
    Extend the radial segments to unit-speed geodesics \(z(t)\) and \(w_i(t)\),
    \(0 \le t \le s / 2\),
    with
    \[
        z(0) = 0, \qquad w_i(0) = v_i, \qquad z(b) = z, \qquad w_i(b) = w_i.
    \]
    The rays \(w_i(t)\) remain in \(Q_i\) since each \(Q_i\) has inradius \(s/2\).
    Similarly,
    \(z(t)\) remains in \(T\).
    Moreover,
    \(w_2(t) = h(w_1(t))\) because \(h\) is an isometry.
    At \(t = s / 2\),
    estimate \eqref{eq:paired-distance-bound} implies that at least one of the distances is at most \(R\).
    Interchanging the indices if necessary,
    we may assume that \(\delta(z(s / 2), w_1(s / 2)) \le R\).
    For \(i = 1, 2\),
    let \(H_i\) and \(K_i\) denote the coefficients \(H\) and \(K\),
    respectively,
    in the distance identity applied to the geodesics \(z(t)\) and \(w_i(t)\).
    Since \(\delta(z(b), w_1(b)) \ge R\), equation \eqref{eq:distance-difference-identity} gives
    \[
        0 \le H_1 \tanh b + K_1 \le H_1 \tanh(s / 2) + K_1 \le 0.
    \]
    As \(b < s / 2\) and \(H_1 \ge 0\), this forces \(H_1 = K_1 = 0\).
    Thus \(z(t)\) moves along \(0 v_1\) toward \(v_1\),
    while \(w_1(t)\) moves away from \(0\) along the same geodesic,
    so \(\delta(0, w_1(s / 2)) = R + s / 2\).
    Since the direction of \(z(t)\) is not along \(0 v_2\),
    we have \(H_2 > 0\).
    Together with \(\delta(z(b), w_2(b)) \ge R\),
    the distance identity yields \(\delta(z(s / 2), w_2(s / 2)) > R\).
    Choose \(z' \in T\) sufficiently close to \(z(s / 2)\) such that \(\delta(0, z') = s / 2\) and \(z'\) does not lie on \(0 v_1\).
    Since \(0\),
    \(z'\) and \(w_1(s / 2)\) are not collinear,
    the strict triangle inequality gives \(\delta(z', w_1(s / 2)) > R\),
    while continuity preserves \(\delta(z', w_2(s / 2)) > R\).
    This contradicts \eqref{eq:paired-distance-bound}.
    We have proved
    \begin{equation}\label{eq:strict-distance-bound}
        0 < b < \frac{s}{2},\ z \in T,\ \delta(0, z) \ge b \qquad \implies \qquad f(z) < R.
    \end{equation}

    Every \(x \in B\) satisfies
    \begin{equation}\label{eq:vertex-distance-bound}
        \min\{d(x, O), d(x, V)\} \le \frac{s}{2},
    \end{equation}
    and equality is possible only at one of the points \(M_j\).
    Indeed,
    the closed hyperbolic disks of radius \(s / 2\) centered at \(0\) and at the vertices \(v_0, \ldots, v_7\) cover \(P\),
    and their boundary circles meet only at the side midpoints.

    Now take any diameter pair \(\{p, q\} \in \calL(B)\),
    and put
    \[
        b = \min\{d(p, O), d(p, V), d(q, O), d(q, V)\} \le \frac{s}{2}.
    \]
    By exchanging the endpoints and using the dual symmetry,
    we may assume that \(b = d(q, V) \le d(p, O)\).
    Choose a lift \(\tilde{p} \in P\) of \(p\).
    Since \(P\) is a Dirichlet domain centered at \(0\),
    we have \(\delta(0, \tilde{p}) = d(O, p)\).
    The rotation of \(P\) through \(\pi / 4\) induces an isometry of \(B\) fixing \(O\) and \(V\).
    Applying a suitable power of this isometry to \(p\) and \(q\) and rotating \(\tilde{p}\) accordingly,
    we may assume that \(\tilde{p} \in T\).
    Choose a lift \(\tilde{q}_1 \in Q_1\) of \(q\) satisfying \(\delta(v_1, \tilde{q}_1) = b\),
    and put \(\tilde{q}_2 \coloneqq h(\tilde{q}_1) \in Q_2\).
    Since \(h\) is a deck transformation,
    \(\tilde{q}_2\) is also a lift of \(q\),
    and \(\delta(v_2, \tilde{q}_2) = b\).
    
    If \(0 < b < s / 2\),
    equation \eqref{eq:strict-distance-bound} contradicts \(d(p, q) = R\).
    
    If \(b = 0\),
    then \(q = V\)
    and \(d(p, V) = R\).
    Choose a lift of \(p\) in a dual Dirichlet octagon centered at a lift of \(V\).
    Since its distance from the center is \(R\),
    it must be a vertex of that octagon.
    All such vertices project to \(O\).
    Hence \(p = O\).
    
    If \(b = s / 2\),
    the definition of \(b\) and equation \eqref{eq:vertex-distance-bound} force both endpoints of the original pair to be among the \(M_j\).
    Let
    \[
        m_j \coloneqq \tanh \frac{s}{4} e^{j \pi i / 4}, \qquad j \in \Z,
    \]
    denote the side midpoints in the disk.
    Since opposite side midpoints are identified,
    \(M_j = \pi(m_j)\),
    with indices taken modulo \(4\).
    The hyperbolic law of cosines gives
    \[
        \cosh \delta(m_j, m_{j + 1}) = \cosh^2 \frac{s}{2} - \sinh^2 \frac{s}{2} \cos \frac{\pi}{4} = \cosh \frac{s}{2}.
    \]
    Hence
    \[
        d(M_j, M_{j + 1}) \le \delta(m_j, m_{j + 1}) = \frac{s}{2} < R.
    \]
    Thus only \(\{O, V\}\),
    \(\{M_0, M_2\}\) and \(\{M_1, M_3\}\) can be diameter pairs.
    
    The six Weierstrass points of \(B\) correspond to the vertices of a regular octahedron in the hyperelliptic quotient \cite[Section 2]{KatzSabourau2006}.
    The rotational symmetries of the octahedron lift to \(\Aut(B)\) and act transitively on the three pairs of opposite vertices.
    Hence
    \[
        \Aut(B) \cdot \{O, V\} = \{\{O, V\}, \{M_0, M_2\}, \{M_1, M_3\}\}.
    \]
    Since \(\{O, V\}\) is a diameter pair by \cite[Proof of Theorem 2]{SBPK} and automorphisms preserve hyperbolic distance,
    all three pairs realize the diameter.
\end{proof}

Let \(P_i = \{p_i, q_i\}\),
\(i = 1, 2, 3\),
denote,
in order,
the three diameter pairs in Proposition \ref{prp:bolza-pairs}.

\begin{lem}\label{lem:bolza-pants}
    There exists a pants decomposition \(\Gamma = \{\gamma_1, \gamma_2, \gamma_3\}\) of \(B\) by simple closed geodesics such that,
    for each \(i = 1, 2, 3\),
    \(\gamma_i\) is the union of two minimizing geodesic segments joining the points of \(P_i\).
\end{lem}

\begin{proof}
    By \cite[Figure 12 (ii)]{SahaSanki},
    the second systoles on \(B\) are explicitly described.
    As illustrated in Figure \ref{fig:bolza-pants},
    we can choose three such geodesics,
    one passing through each diameter pair \(P_i\),
    so that they are pairwise disjoint.
    \begin{figure}
        \centering
        \resizebox{0.4\linewidth}{!}{\begin{tikzpicture}[scale=3.0, line cap=round, line join=round]
  \definecolor{bolzared}{RGB}{214,39,39}
  \definecolor{bolzablue}{RGB}{24,71,211}
  \definecolor{bolzagreen}{RGB}{19,148,79}
  
  \pgfmathsetmacro{\rhoV}{sqrt((1+sqrt(2))/(2+sqrt(2)))}
  \pgfmathsetmacro{\rhoM}{sqrt(sqrt(2)-1)}
  \pgfmathsetmacro{\halfDiag}{\rhoM/sqrt(2)}
  \pgfmathsetmacro{\sideCenter}{(1+\rhoV*\rhoV)/(2*\rhoV*cos(22.5))}
  \pgfmathsetmacro{\pairT}{2*\rhoM/(1+\rhoM*\rhoM)}
  \pgfmathsetmacro{\edgeSlope}{(\sideCenter-\pairT)/\rhoM}
  \pgfmathsetmacro{\edgeHeight}{(\sideCenter-sqrt(\sideCenter*\sideCenter-1-\edgeSlope*\edgeSlope))/(1+\edgeSlope*\edgeSlope)}
  \pgfmathsetmacro{\edgeOffset}{\edgeSlope*\edgeHeight}

  \newcommand{\HyperbolicArc}[5]{%
    \begingroup
      \pgfmathsetmacro{\xa}{#1}%
      \pgfmathsetmacro{\ya}{#2}%
      \pgfmathsetmacro{\xb}{#3}%
      \pgfmathsetmacro{\yb}{#4}%
      \pgfmathsetmacro{\detA}{2*(\xa*\yb-\xb*\ya)}%
      \pgfmathsetmacro{\qa}{1+\xa*\xa+\ya*\ya}%
      \pgfmathsetmacro{\qb}{1+\xb*\xb+\yb*\yb}%
      \pgfmathsetmacro{\cx}{(\qa*\yb-\qb*\ya)/\detA}%
      \pgfmathsetmacro{\cy}{(\xa*\qb-\xb*\qa)/\detA}%
      \pgfmathsetmacro{\cr}{sqrt(\cx*\cx+\cy*\cy-1)}%
      \pgfmathsetmacro{\angA}{atan2(\ya-\cy,\xa-\cx)}%
      \pgfmathsetmacro{\angB}{atan2(\yb-\cy,\xb-\cx)}%
      \pgfmathsetmacro{\angleDelta}{mod(\angB-\angA+540,360)-180}%
      \draw[#5] (\xa,\ya) arc[start angle=\angA, delta angle=\angleDelta, radius=\cr];%
    \endgroup
  }

  \foreach \k in {0,...,7} {
    \pgfmathsetmacro{\a}{22.5+45*\k}
    \pgfmathsetmacro{\b}{22.5+45*(\k+1)}
    \HyperbolicArc{\rhoV*cos(\a)}{\rhoV*sin(\a)}%
                  {\rhoV*cos(\b)}{\rhoV*sin(\b)}%
                  {black,line width=1.15pt}
  }

  \draw[bolzagreen, line width=1.55pt]
    ({\rhoV*cos(202.5)},{\rhoV*sin(202.5)}) --
    ({\rhoV*cos(22.5)},{\rhoV*sin(22.5)});

  \HyperbolicArc{-\rhoM}{0}{0}{\rhoM}{bolzared,line width=1.55pt}
  \HyperbolicArc{0}{-\rhoM}{\rhoM}{0}{bolzared,line width=1.55pt}

  \HyperbolicArc{-\halfDiag}{\halfDiag}{-\edgeOffset}{\edgeHeight}%
                {bolzablue,line width=1.55pt}
  \HyperbolicArc{\edgeOffset}{\edgeHeight}{\halfDiag}{\halfDiag}%
                {bolzablue,line width=1.55pt}
  \HyperbolicArc{-\halfDiag}{-\halfDiag}{-\edgeOffset}{-\edgeHeight}%
                {bolzablue,line width=1.55pt}
  \HyperbolicArc{\edgeOffset}{-\edgeHeight}{\halfDiag}{-\halfDiag}%
                {bolzablue,line width=1.55pt}

  \tikzset{bolza dot/.style={circle,draw=black,fill=bolzared,
                             line width=0.2pt,inner sep=0pt,minimum size=2.4pt}}
  \foreach \x/\y in {
     0/0,
     0/\rhoM, 0/-\rhoM,
     \rhoM/0, -\rhoM/0,
     \halfDiag/\halfDiag, -\halfDiag/\halfDiag,
     \halfDiag/-\halfDiag, -\halfDiag/-\halfDiag,
     \edgeOffset/\edgeHeight, -\edgeOffset/\edgeHeight,
     \edgeOffset/-\edgeHeight, -\edgeOffset/-\edgeHeight
  } {
     \node[bolza dot] at (\x,\y) {};
  }
  \node[bolza dot] at ({\rhoV*cos(22.5)},{\rhoV*sin(22.5)}) {};
  \node[bolza dot] at ({\rhoV*cos(202.5)},{\rhoV*sin(202.5)}) {};

  \node[font=\scriptsize,above right,inner sep=1pt] at (0,0) {$O$};
  \node[font=\scriptsize,above left,inner sep=1.6pt] at (0,\rhoM) {$M_2$};
  \node[font=\scriptsize,right,inner sep=1.6pt] at (\halfDiag,\halfDiag) {$M_1$};
  \node[font=\scriptsize,left,inner sep=1.6pt] at (-\halfDiag,\halfDiag) {$M_3$};
  \node[font=\scriptsize,right,inner sep=1.6pt] at (\rhoM,0) {$M_0$};
  \node[font=\scriptsize,right,inner sep=1.6pt]
       at ({\rhoV*cos(22.5)},{\rhoV*sin(22.5)}) {$V$};
\end{tikzpicture}}
        \caption{Three pairwise disjoint simple closed geodesics on the Bolza surface.
                 The green,
                 red and blue curves correspond to \(P_1\),
                 \(P_2\) and \(P_3\),
                 respectively.}
        \label{fig:bolza-pants}
    \end{figure}
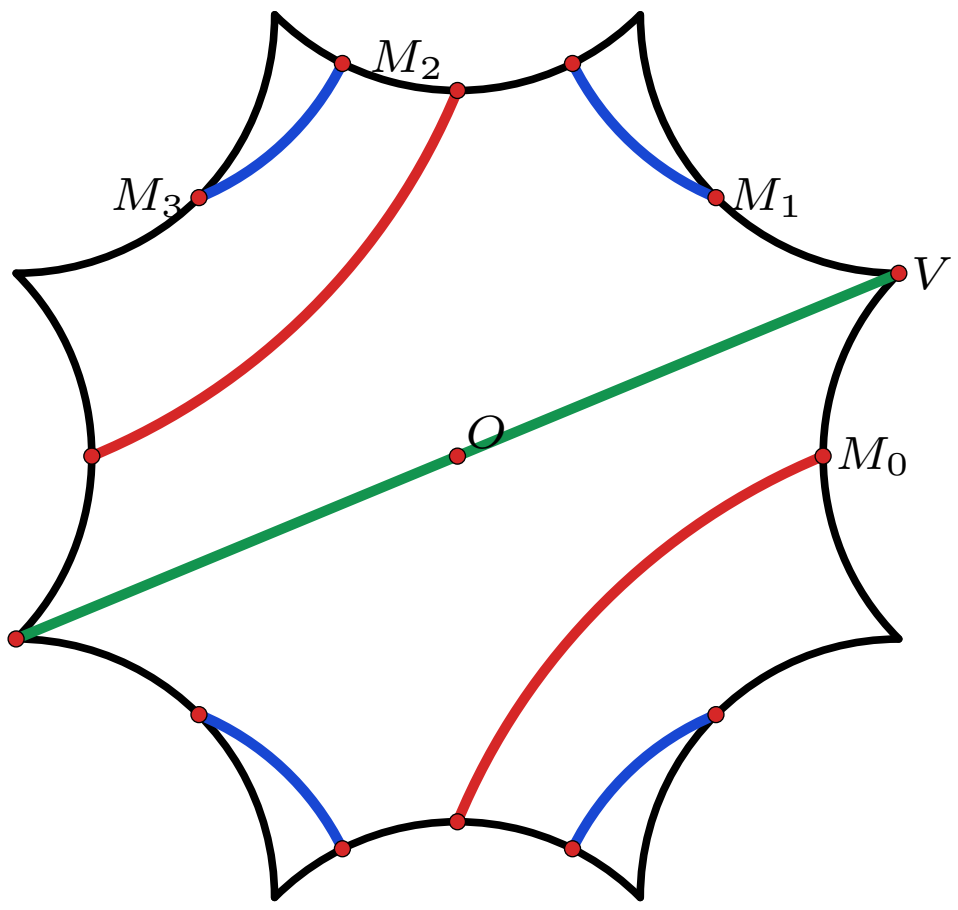
    By \cite[Proof of Theorem 1.3]{SahaSanki},
    each \(\gamma_i\) has length \(2 D_B\).
    Since \(\gamma_i\) passes through \(p_i\) and \(q_i\),
    and \(d_B(p_i, q_i) = D_B\),
    the two subarcs joining these points must each have length \(D_B\).
    Hence both subarcs are minimizing geodesic segments.
    Therefore the three pairwise disjoint simple closed geodesics form a pants decomposition of \(B\).
\end{proof}

\begin{proof}[Proof of Corollary \ref{cor:bolza}]
    By Proposition \ref{prp:bolza-pairs},
    \(\calL(B) = \{P_1, P_2, P_3\}\).
    Lemma \ref{lem:bolza-pants} assigns to these pairs three disjoint simple closed geodesics,
    each formed by two minimizing segments joining the corresponding pair.
    Theorem \ref{thm:main} now gives the result.
\end{proof}

\enlargethispage{2\baselineskip}
{\small
\bibliographystyle{amsplain}
\bibliography{bolza_nonmini_references}
}
\end{document}